\documentstyle[12pt]{article}
\newcommand{\lap}{\mbox{$\bigtriangleup$}}

\newtheorem{mthm}{Theorem}

\begin{document}
\title
{Uniqueness of positive bound states to Schrodinger systems with
critical exponents}
\author{ Congming Li
\thanks{Partially supported by NSF Grant DMS-0401174 }
\hspace{.2in} Li Ma\thanks{Partially supported by the National
Natural Science Foundation of China 10631020 and SRFDP 20060003002 }
}
\date{}
\maketitle
\bigskip

\begin{abstract}
We prove the uniqueness for the positive solutions of the following
elliptic systems:
\begin{eqnarray*}
\left\{\begin{array}{ll} - \lap (u(x) ) = u(x)^{\alpha}v(x)^{\beta}
 \\
- \lap ( v(x) ) = u(x)^{\beta} v(x)^{\alpha}
\end{array}
\right.
\end{eqnarray*}
Here $x\in R^n$, $n\geq 3$, and $1\leq \alpha, \beta\leq
\frac{n+2}{n-2}$ with $\alpha+\beta=\frac{n+2}{n-2}$. In the special
case when $n=3$ and $\alpha =2, \beta=3$, the systems come from the
stationary Schrodinger system with critical exponents for
Bose-Einstein condensate. As a key step, we prove the radial
symmetry of the positive solutions to the elliptic system above with
critical exponents.
\end{abstract}


\bigskip

\newpage

\section{Introduction}

In this paper, we consider the uniqueness of positive solutions to
the following stationary Schrodinger system:
\begin{equation} \label{EllS}
\left\{\begin{array}{ll} - \lap (u(x) ) = u(x)^{\alpha}v(x)^{\beta}
 \\
- \lap ( v(x) ) = u(x)^{\beta} v(x)^{\alpha}
\end{array}
\right.
\end{equation}
Here $x\in R^n$, $n\geq 3$, and $1\leq \alpha, \beta\leq
\frac{n+2}{n-2}$ with $\alpha+\beta=\frac{n+2}{n-2}$. In the special
case when $n=3$ and $\alpha =2, \beta=3$, the systems come from the
stationary Schrodinger system with critical exponents for
Bose-Einstein condensate (\cite{KL},\cite{LW1},\cite{LW2}, and
\cite{MZ}). In the earlier works \cite{LW1},\cite{LW2}, and
\cite{MZ}, people pay more attention to the elliptic system
(\ref{EllS}) with subcritical exponents. Very interestingly, Chen
and Li have proved that the best constant in weighted
Hardy-Littlewood-Sobolev inequality can be achieved by explicit
radially symmetric functions (see \cite{CL3} and \cite{L}). As a
consequence of their work, the uniqueness of positive solutions to
the corresponding elliptic system (it is (\ref{EllS}) in the case
when $\alpha=0$ and $\beta=\frac{n+2}{n-2})$) has been settled down.
However, when $0<\alpha, \beta$, the uniqueness of smooth positive
solutions to the stationary Schrodinger system (\ref{EllS}) is an
open question. Generally speaking, there are very few result even
for the uniqueness of positive solutions to the ordinary
differential systems. The aim of this paper is to prove the radial
symmetry and uniqueness of positive solutions to (\ref{EllS}) with
critical exponents and $1\leq \alpha<\beta\leq \frac{n+2}{n-2}$.

As one can expect, just like in the work M.Weinstein \cite{We1} in
the scalar case with the sub-critical exponent, that there is a
closed relationship between the stationary Schrodinger system with
critical exponent with the Hardy-Littlewood-Sobolev inequality. As
we show below, this is true.

 Since we shall use Hardy-Littlewood-Sobolev
 inequality to prove radial symmetry of our solutions,
 let's do an excursion about recent progress
 of Lieb's conjecture. Let us begin by recalling the well-known Hardy-Littlewood-Sobolev inequalities.
 Let $0 <\lambda < n$, $1 < s, r < \infty$,
and $\|f\|_p$ be the $L^p(R^n)$ norm of the function $f$. We shall
write by  $\|f\|_{p,\Omega}$ the $L^p$ norm of the function on the
domain $\Omega$. Then, classical Hardy-Littlewood-Sobolev inequality
(HLS) states that:
\begin{equation}
    \int_{R^n}\int_{R^n} \frac{f(x) g(y)}{|x - y|^{\lambda}} dx dy \leq
    C_{s,\lambda, n} \|f\|_r\|g\|_s
    \label{HLS}
\end{equation} for any $f \in L^r(R^n)$, $g \in L^s(R^n)$, and for
$\frac{1}{r} + \frac{1}{s} + \frac{\lambda}{n} = 2$. Hardy and
Littlewood also introduced the double weighted inequality, which was
later generalized by Stein and Weiss  in \cite{SW2} in the following
form:
\begin{equation}
   \left|\int_{R^n}\int_{R^n}
   \frac{f(x)g(y)}{|x|^{\alpha_0}|x-y|^{\lambda}|y|^{\beta_0}}dxdy
   \right| \le C_{\alpha_0,\beta_0,s,\lambda,n }\|f\|_r \|g\|_s
   \label{whls}
\end{equation}
where $\alpha_0+\beta_0 \geq 0$,
\begin{equation}
1 - \frac{1}{r} - \frac{\lambda}{n} < \frac{\alpha_0}{n} < 1 -
\frac{1}{r} , \; \mbox{ and } \frac{1}{r} + \frac{1}{s} +
\frac{\lambda + \alpha_0 + \beta_0}{n} = 2 . \label{B1}
\end{equation}

The best constant in the weighted inequality (\ref{whls}) can be
obtained by maximizing the functional
\begin{equation}
    J(f,g) =
    \int_{R^n}\int_{R^n}\frac{f(x)g(y)}
    {|x|^{\alpha_0}|x-y|^{\lambda}|y|^{\beta_0}} dx dy
\label{J}
\end{equation}
under the constraints $ \|f\|_r = \|g\|_s =1 $. Then the
corresponding Euler-Lagrange equations are the system of integral
equations:
\begin{equation}
 \left \{
   \begin{array}{l}
      \lambda_1 r {f(x)}^{r-1} = \frac{1}{|x|^{\alpha_0}}\int_{R^{n}} \frac{g(y)}{|y|^{\beta_0}|x-y|^{\lambda}}  dy\\
      \lambda_2 s {g(x)}^{s-1} = \frac{1}{|x|^{\beta_0}}\int_{R^{n}} \frac{f(y)}{|y|^{\alpha_0}|x-y|^{\lambda}}  dy
   \end{array}
   \right.
   \label{DHLSEL}
\end{equation} where $f, g \geq 0, \;x \in R^n$,
and $\lambda_1 r =\lambda_2 s =J(f,g)$.

Let $u=c_1 f^{r-1}$, $v=c_2 g^{s-1}$, $p=\frac{1}{r-1}$,
$q=\frac{1}{s-1}$, $pq \neq 1 $, and for a proper choice of
constants $c_1$ and $c_2$, system (\ref{DHLSEL}) becomes
\begin{equation}
 \left \{
   \begin{array}{l}
      u(x) = \frac{1}{|x|^{\alpha}}\int_{R^{n}} \frac{v(y)^q}{|y|^{\beta_0}|x-y|^{\lambda}}  dy\\
      v(x) = \frac{1}{|x|^{\beta}}\int_{R^{n}} \frac{u(y)^p}{|y|^{\alpha_0}|x-y|^{\lambda}}  dy
   \end{array}
   \right.\label{ELWHLS}
 \end{equation} where $u, v \geq 0$, $0 <p, q<\infty$, $0 <\lambda < n$,  $\frac{\alpha_0}{n}  <
  \frac{1}{p+1}<\frac{\lambda+\alpha_0}{n}$, and
 $\frac{1}{p+1}+\frac{1}{q+1}=\frac{\lambda+\alpha_0+\beta_0}{n}$.

Note that in the special case where $\alpha_0 =0$ and $\beta_0 =0$,
system (\ref{ELWHLS}) reduces to the following system:

\begin{equation}
 \left \{
 \begin{array}{l}
      u(x) = \int_{R^{n}} \frac{v^q(y)}{|x - y|^{\lambda}} dy\\
      v(x) = \int_{R^{n}} \frac{u^p(y)}{|x - y|^{\lambda}} dy
\end{array}
\right. \label{sys} \end{equation} with \begin{equation}
\frac{1}{q+1}+\frac{1}{p+1}=\frac{\lambda}{n}. \label{a1}
\end{equation}

It is well-known that this integral system is closely related to the
system of partial differential equations

\begin{equation}
 \left \{
 \begin{array}{l}
      (-\Delta)^{\gamma/2} u = v^q , \;\; u>0, \mbox{ in } R^n,\\
      (-\Delta)^{\gamma/2} v = u^p, \;\; v>0, \mbox{ in } R^n,
\end{array}
\right. \label{des} \end{equation} where $\gamma = n - \lambda$.

When $p=q=\frac{n+\gamma}{n-\gamma}$, and $u(x)=v(x)$, system
(\ref{sys}) becomes the single equation:

\begin{equation}
      u(x) = \int_{R^{n}} \frac{u(y)^{\frac{n+\gamma}{n-\gamma} }}{|x - y|^{n - \gamma}}
      dy,\;\;
       u >0,  \mbox{ in } R^n.
\label{eq} \end{equation}

The corresponding PDE is the well-known family of semi-linear
equations

\begin{equation} (- \Delta)^{\gamma/2} u = u^{(n+\gamma)/(n-\gamma)} , \;\; u>0,
\;\; \mbox{ in } R^n \label{de} \end{equation}

In particular, when $n \geq 3 $, and $\gamma = 2,$ (\ref{de})
becomes

\begin{equation} -\Delta u =  u^{ (n+2)/(n-2) }, \;\; u>0, \mbox{ in } R^n.
\label{1.2} \end{equation}

The classification of the solutions of (\ref{1.2}) has provided an
important ingredient in the study of the well-known Yamabe problem
and the prescribing scalar curvature problem. Equation (\ref{1.2})
were studied by Gidas, Ni, and Nirenberg \cite{GNN}, Caffarelli,
Gidas, and Spruck \cite{CGS},  Chen and Li \cite{CL}, and Li
\cite{Li}. They classified all the solutions. Recently, Wei and Xu
\cite{WX} generalized this result to the solutions of more general
equation (\ref{de}) with $\gamma$ being any even numbers between 0
and n.

Although the systems for other real values of $\alpha, \beta$
between 0 and n are of interest to people, we shall only concentrate
in this paper to the system (\ref{EllS}) with critical exponents
when $1\leq \alpha, \beta\leq \frac{n+2}{n-2}$ and $\alpha+\beta=
\frac{n+2}{n-2}$.

Our main results are

\begin{mthm}\label{thm1}. Any $L^{\frac{2n}{n-2}}(\mathbf{R}^n)\times
L^{\frac{2n}{n-2}}(\mathbf{R}^n)$ positive solution pair $(u,v)$ to
the system (\ref{EllS}) with critical exponents are radial symmetric
functions.
\end{mthm}
and

\begin{mthm}\label{thm2}.
Assume that $1\leq \alpha<\beta\leq \frac{n+2}{n-2}$. Then any
$L^{\frac{2n}{n-2}}(\mathbf{R}^n)\times
L^{\frac{2n}{n-2}}(\mathbf{R}^n)$ radial symmetric solution pair
$(u,v)$ to the system (\ref{EllS}) with critical exponents are
unique such that $u=v$.
\end{mthm}

We point out that when $u=v$, the elliptic system (\ref{EllS})
reduces to the elliptic equation with critical exponent (\ref{1.2}).
Then $u=v$ is in a special family of functions:
\begin{equation}
\phi_{x_o,t}(x)=c(\frac{t}{t^2 + |x - x_o|^2})^{(n-2)/2}
\label{solu}
\end{equation}
where $t>0$, $x_o\in\mathbf{R}^n$, with some positive constants c
such that each $\phi_{x_o,t}(x)$ solves (\ref{1.2}). This family of
functions are important in the study of (\ref{EllS}).

Our results are motivated from the previous work \cite{CLO2},where
Chen, Li, and Ou considered more general system (\ref{sys}) and
established the symmetry and monotonicity of the solutions. In
\cite{CL2}, Chen and Li also obtained a regularity result of the
solutions to (\ref{sys}). To establish the symmetry of the solution
to (\ref{sys}), Chen, Li, and Ou \cite{CLO} \cite{CLO2} \cite{CLO4}
introduced a new idea, an integral form of the method of moving
planes. It is entirely different from the traditional method used
for partial differential equations. Instead of relying on maximum
principles, certain integral norms were estimated. The new method is
a very powerful tool in studying qualitative properties of other
integral equations and systems. In fact, following Chen, Li, and
Ou's work, Jin and Li \cite{CJ} studied the symmetry of the
solutions to the more general system (\ref{ELWHLS}).

Chen and Ma \cite{MC} discussed the Liouville type theorem for the
positive solutions to the elliptic system (\ref{des}).

In this paper, we first prove the radial symmetry of the solutions
to (\ref{EllS}) with critical exponents. It is obvious that the
radial symmetry of the solutions reduces (\ref{EllS}) to a system of
ODEs, which has special solution pair
$(\phi_{o,t}(x),\phi_{o,t}(x))$. To prove the uniqueness, we prove
that $u(0)=v(0)$. Then by the uniqueness of the initial value
problem for ODE, we conclude that $u=v=\phi_{o,t}$. This is the key
observation in establishing the uniqueness of positive solutions for
(\ref{EllS}) with critical exponents.

Theorems \ref{thm1} and \ref{thm2} will be proved in the next two
sections.

\section{Proof of the Radial symmetry}
We use the moving plane method introduced by Chen-Li-Ou in
\cite{CLO}. We shall use the Hardy-Littlewood-Sobolev inequality:
\begin{equation} \label{HLS3}
|Tf|_p\leq |f|_{\frac{np}{n+2p}},
\end{equation}
where $C(n,p)$ is a uniform positive constant and
$$
Tf(x)=\int_{\mathbf{R}^n}|x-y|^{2-n}f(y)dy.
$$

{\bf The Proof of Theorem \ref{thm1}}. For each $\lambda\in
\mathbf{R}$, we denote by
$$
H_{\lambda}=\{x\in\mathbf{R}^n; x_1<\lambda\}.
$$

For each $x=(x_1,x')\mathbf{R}^n$, we let
$$
x_{\lambda}=(2\lambda-x_1,x')
$$
be the reflection point of $x$ with respect to the hyperplane
$\partial H_{\lambda}$. We let $e_1=(1,0,...0)$.

We define
$$
u_{\lambda}(x)=u(x_{\lambda}), \; \; B_{\lambda}^u=\{x\in
H_{\lambda}; u_{\lambda}(x)>u(x)\},
$$
and
$$
v_{\lambda}(x)=v(x_{\lambda}), \;\; B_{\lambda}^v=\{x\in
H_{\lambda}; v_{\lambda}(x)>v(x)\}.
$$

To do the moving plane method, we need the following formula, which
is obtained by a change of variables.
$$
u(x)=\int_{H_{\lambda}}\frac{u^{\alpha}v^{\beta}(y)}{|x-y|^{n-2}}dy+
\int_{H_{\lambda}}
\frac{u_{\lambda}^{\alpha}v_{\lambda}^{\beta}(y)}{|x_{\lambda}-y|^{n-2}}dy,
$$
and
$$
u_{\lambda}(x)=\int_{H_{\lambda}}\frac{u^{\alpha}v^{\beta}(y)}{|x_{\lambda}-y|^{n-2}}dy+
\int_{H_{\lambda}}\frac{u_{\lambda}^{\alpha}v_{\lambda}^{\beta}(y)}{|x-y|^{n-2}}dy.
$$
Then we have
\begin{equation} \label{equ3}
u_{\lambda}(x)-u(x)=\int_{H_{\lambda}}
(u_{\lambda}^{\alpha}v_{\lambda}^{\beta}-u^{\alpha}v^{\beta})(y)
(\frac{1}{|x-y|^{n-2}}-\frac{1}{|x_{\lambda}-y|^{n-2}})dy.
\end{equation}
Note that for $x\in H_{\lambda}$, we have
$$
\frac{1}{|x-y|^{n-2}}>\frac{1}{|x_{\lambda}-y|^{n-2}}.
$$
Then for $x\in B^u_{\lambda}$, we have
\begin{equation} \label{est1}
\left\{\begin{array}{ll}
0&\leq  u_{\lambda}(x)-u(x)\\
&\leq
 \alpha\int_{B^u_{\lambda}}
 \frac{u_{\lambda}^{\alpha-1}v_{\lambda}^{\beta}(u_{\lambda}-u)}{|x-y|^{n-2}}dy
 +\beta
 \int_{B^v_{\lambda}}
 \frac{u_{\lambda}^{\alpha}v_{\lambda}^{\beta-1}(v_{\lambda}-v)}{|x-y|^{n-2}}dy\\
 &:=I+II
 \end{array}
 \right.
 \end{equation}

Let $p=\frac{2n}{n-2}$. Using the Hardy-Littlewood-Sobolev
inequality (\ref{HLS3}) we can bound the first term $I$ in
(\ref{est1}) by
\begin{equation} \label{est2}
\left\{\begin{array}{ll}  |I|_p\leq
C(n,p)|u_{\lambda}^{\alpha-1}v_{\lambda}^{\beta}(u_{\lambda}-u)|_{\frac{2n}{n+2}}\\
\leq
C(n,p)|u_{\lambda}|_p^{\alpha-1}|v_{\lambda}|_p^{\beta}|u_{\lambda}-u|_p
\end{array}
 \right.
 \end{equation}
Here the integrations are over the set $B^u_{\lambda}$.

Using again the Hardy-Littlewood-Sobolev inequality (\ref{HLS3}) we
can bound the first term $II$ in (\ref{est1}) by
\begin{equation} \label{est3}
\left\{\begin{array}{ll}  |II|_p\leq
C(n,p)|u_{\lambda}^{\alpha}v_{\lambda}^{\beta-1}(v_{\lambda}-v)|_{\frac{2n}{n+2}}\\
\leq
C(n,p)|u_{\lambda}|_p^{\alpha}|v_{\lambda}|_p^{\beta-1}|v_{\lambda}-v|_p
\end{array}
 \right.
 \end{equation}
Here the integrations are over the domain $B^v_{\lambda}$. Hence, we
have
\begin{equation}\label{mov1}
|u_{\lambda}-u|_{p, B^u_{\lambda}}\leq
C(n,p)(|u_{\lambda}|_{p,B^u_{\lambda}}^{\alpha-1}|v_{\lambda}|_{p,B^u_{\lambda}}^{\beta}|u_{\lambda}-u|_{p,B^u_{\lambda}}
+|u_{\lambda}|_{p,B^v_{\lambda}}^{\alpha}|v_{\lambda}|_{p,B^v_{\lambda}}^{\beta-1}|v_{\lambda}-v|_{p,B^v_{\lambda}})
\end{equation}

 Similarly, we have for following formulae for $v$ and
$v_{\lambda}$.
$$
v(x)=\int_{H_{\lambda}}\frac{u^{\alpha}v^{\beta}(y)}{|x-y|^{n-2}}dy+
\int_{H_{\lambda}}
\frac{v_{\lambda}^{\alpha}u_{\lambda}^{\beta}(y)}{|x_{\lambda}-y|^{n-2}}dy,
$$
and
$$
v_{\lambda}(x)=\int_{H_{\lambda}}\frac{v^{\alpha}u^{\beta}(y)}{|x_{\lambda}-y|^{n-2}}dy+
\int_{H_{\lambda}}\frac{v_{\lambda}^{\alpha}u_{\lambda}^{\beta}(y)}{|x-y|^{n-2}}dy.
$$
Then we have the following estimate
\begin{equation}\label{mov2}
|v_{\lambda}-v|_p\leq
C(n,p)(|v_{\lambda}|_{p,B^v_{\lambda}}^{\alpha-1}|u_{\lambda}|_{p,B^v_{\lambda}}^{\beta}|v_{\lambda}-v|_{p,B^v_{\lambda}}
+|v_{\lambda}|_{p,B^u_{\lambda}}^{\alpha}|u_{\lambda}|_{p,B^u_{\lambda}}^{\beta-1}|u_{\lambda}-u|_{p,B^u_{\lambda}})
\end{equation}

  After these preparations, we can use the moving plane method as
  developed in \cite{CLO} to prove the radial symmetry of the solutions.

  At first, let's start the plane from the infinity.
  Indeed, for $\lambda>>1$ large enough, we know that the quantities
  $$
|v_{\lambda}|_{p,B^u_{\lambda}}, |u_{\lambda}|_{p,B^u_{\lambda}},
|v_{\lambda}|_{p,B^v_{\lambda}}
  $$
  and $$
|u_{\lambda}|_{p,B^v_{\lambda}}
  $$
  all are small, which give us that
  $$
|u_{\lambda}-u|_{p, B^u_{\lambda}} \leq
\frac{1}{2}|v_{\lambda}-v|_{p,B^v_{\lambda}}
$$
and
$$
|v_{\lambda}-v|_{p,B^v_{\lambda}}\leq
\frac{1}{2}|u_{\lambda}-u|_{p,B^u_{\lambda}}.
$$
These imply that $|u_{\lambda}-u|_{p, B^u_{\lambda}} =0$ and
$|v_{\lambda}-v|_{p,B^v_{\lambda}}=0$. These say that
$B^u_{\lambda}=\phi$ and $B^v_{\lambda}=\phi$.

Next we define
$$
\lambda_0=\{\lambda\in\mathbf{R}; B^u_{\lambda'}=\phi\; \mbox{for
all}\; \lambda'\geq\lambda\}.
$$
Then it follows from the fact that $u(x)\to 0$ as $|x|\to \infty$
and $u(x)>0$ in $\mathbf{R}^n$ that $\lambda<+\infty$. By the
definition of $\lambda_0$, we have $u_{\lambda_0}(x)\leq u(x)$ foe
$x\in H_{\lambda_0}$. Using the expression (\ref{equ3}), we see that
$u_{\lambda_0}(x)<u(x)$ for $x\in H_{\lambda_0}$. This implies that
$|2\lambda e_1-B^{u}_{\lambda}|\to 0$ as $\lambda\to \lambda_0$. It
is now standard to know (see \cite{CLO}) that $u_{\lambda_0}=u$,
 which then gives us $v_{\lambda_0}=v $. Since $x_1$ can be any
 directions, we conclude that $u$ and $v$ are radial symmetric about
 some point $x_o$.

\section{Proof of the Uniqueness}

In some sense, the proof of Theorem \ref{thm2} is just at hand by
using the integral expression of the solution pair $(u,v)$.

{\bf Proof of Theorem \ref{thm2}}. Let $(u, v)\in
L^{\frac{2n}{n-2}}(\mathbf{R}^n)\times
L^{\frac{2n}{n-2}}(\mathbf{R}^n)$ be a pair of solutions to system
(\ref{EllS}). By Theorem \ref{thm1}, we know that $u$ and $v$ are
radial symmetric about the some point $x_0$. We may say, $x_0=0$.

Recall that we have assumed that
$1\leq\alpha<\beta<\frac{n+2}{n-2}$. Since, $u\in
L^{\frac{2n}{n-2}}(\mathbf{R}^n)$ and $v\in
L^{\frac{2n}{n-2}}(\mathbf{R}^n)$. Using the same method in
\cite{CLO}, we have that $u\in C^2(\mathbf{R}^n)$ and  $v\in
C^2(\mathbf{R}^n)$ with
$$
u(x)\to 0, v(x)\to 0,
$$
as $|x|\to \infty$.

Since our solution $u$ is radially symmetric, hence we can write, in
polar coordinates, the first equation in (\ref{EllS}) as
$$ (r^{n-1} u'(r))' = - r^{n-1} u(r)^{\alpha}v(r)^{\beta} ,$$
where $r = |x|$.

Integrating both sides from $0$ to $r$ yields
$$ r^{n-1} u'(r) = - \int_0^r s^{n-1} u^{\alpha}v^{\beta} (s) ds .$$
It follows by another integration that
\begin{equation} u(r) =
u(0) - \int_0^r \frac{1}{\tau^{n-1}} \int_0^{\tau}
s^{n-1}u^{\alpha}v^{\beta} ds d\tau . \label{7}
\end{equation}

Similarly, for $v(r)$, we have
\begin{equation}
v(r) = v(0) - \int_0^r \frac{1}{\tau^{n-1}} \int_0^{\tau} s^{n-1}
v^{\alpha}u^{\beta} ds d\tau . \label{8}
\end{equation}

As we mentioned in the introduction, we need only to show that $u(0)
= v(0)$. Otherwise, suppose
\begin{equation}
u(0) < v(0) , \label{9}
\end{equation}
then by continuity, for all small $r>0$,
\begin{equation}
u(r) < v(r) . \label{10}
\end{equation}
In other word, there exists an $R > 0$, such that
\begin{equation}
u(r) < v(r) , \;\; \forall r \in (0, R). \label{11}
\end{equation}

Let $R_o$ be the supreme value of $R$, such that (\ref{11}) holds.
Then $R_o\leq \infty$ and $u(R_o)=v(R_o)$, where we have used the
fact that $u(+\infty)=v(+\infty)=0$. By the definition of $R_o$ and
$\alpha<\beta$, we have that
\begin{equation}
u(r)^{\alpha}v(r)^{\beta} > v(r)^{\alpha}u(r)^{\beta} , \;\; \forall
r \in (0, R_o). \label{12}
\end{equation}

Then we have from (\ref{7}) and (\ref{8}) that
$$0>u(0)-v(0)=
\int_0^{R_o} \frac{1}{\tau^{n-1}} \int_0^{\tau} s^{n-1}
(u^{\alpha}v^{\beta}-u^{\beta}v^{\alpha})(s) ds d\tau>0.
$$
 This is impossible.

 Similarly, one can
show that $u (0) > v(0)$ is impossible. Therefore, we must have
$$ u(0) = v(0). $$

Finally, by the standard ODE theory, we arrive at
$$ u(r)\equiv v(r) . $$
Hence, our elliptic system (\ref{EllS}) has been reduced to the
elliptic equation with critical exponent (\ref{1.2}). By now, it is
standard to know that our solutions pair $u$ and $v$ are of the form
(\ref{solu}). This completes the proof of Theorem \ref{thm2}.

{\em Authors' Addresses and E-mails:}
\medskip

Congming Li

Department of Applied Mathematics

Campus Box 526, University of Colorado at Boulder,

Boulder CO 80309,USA

cli@colorado.edu

\medskip

Li Ma

Department of Mathematics

Tsinghua University

Beijing 100084

China

lma@math.tsinghua.edu.cn


\bigskip



\begin{thebibliography}{CL}

\bibitem{CGS} L. Caffarelli, B. Gidas, and J. Spruck, \,\,
\textit{Asymptotic symmetry and local behavior of semilinear
elliptic equations with critical Sobolev growth}, \,\, Comm. Pure
Appl. Math. XLII,  (1989), 271--297

\bibitem{CL} W. Chen and C. Li, \,\, \textit{Classification of solutions of some
nonlinear elliptic equations}, \,\, Duke Math. J.,  {\bf 63}
(1991), 615--622.

\bibitem{CL2} W. Chen and C. Li, \,\, \textit{Regularity of
Solutions for a system of Integral Equations}, Comm. Pure and
Appl. Anal., 4(2005), 1--8.

\bibitem{CL3} W. Chen and C. Li, \,\,
\textit{The best constant in weighted Hardy-Littlewood-Sobolev
inequality}, Proc. AMS, accepted, 2007.

\bibitem{CLO} W. Chen, C. Li, and B. Ou, \,\, \textit{Classification
of solutions for an integral equation}, \,\,  Comm. Pure and Appl.
Math., 59(2006), 330--343.

\bibitem{CLO2} W. Chen, C. Li, and B. Ou, \,\, \textit{Classification
of solutions for a system of integral equations}, \,\, Comm. in
Partial Differential Equations, 30(2005), 59--65


\bibitem{CLO4} W. Chen, C. Li, and B. Ou, \,\, \textit{Qualitative Properties
of Solutions for an Integral Equation}, \,\, Disc. \& Cont. Dynamics
Sys., 12(2005), 347--354.

\bibitem{GNN} B. Gidas, W. M. Ni, and L. Nirenberg, \,\, \textit{Symmetry of positive solutions
of nonlinear elliptic equations in $R^{n},$ } \,\, Mathematical
Analysis and Applications,  Academic Press, New York, 1981.

\bibitem{CJ}  C. Jin and C. Li, \,\, \textit{Symmetry of Solutions to Some Integral
Equations}, Proc. Amer. Math. Soc., 134(2006), 1661-1670.

\bibitem{CJ2}  C. Jin and C. Li, \,\, \textit{Quantitative Analysis of Some System of Integral Equations},
Cal. Var. PDEs, 26(2006), p447-457.

\bibitem{KL}
Kanna, T. and Lakshmanan, M., {\em Exact soliton solutions, shape
changing collisions, and partially coherent solitons in coupled
nonlinear Schr\"{o}dinger equations}, Phys. Rev. Lett.,
86(2001)5043.


\bibitem{Li} C. Li, \,\, \textit{Local asymptotic symmetry of singular
solutions to nonlinear elliptic equations}, \,\, Invent. Math.,
 {\bf 123} (1996), 221--231.

\bibitem{L} E. Lieb, \,\, \textit{Sharp constants in the Hardy-Littlewood-Sobolev and
related inequalities}, \,\, Ann. of Math., {\bf 118}(1983),
349--374.

\bibitem{LLim} C. Li and J. Lim, \,\, \emph{The singularity analysis of
solutions to some integral edquations}, CPAA, in press.

\bibitem{LW1}
Lin, T.C. and Wei, J., {\em Ground state of $N$ coupled nonlinear
Schr\"{o}dinger equations in $\mathbf{R}^n$, $n\leq 3$}, Commun.
Math. Phys., 255(2005)629-653.


\bibitem{LW2}
Lin, T.C. and Wei, J., {\em Spikes in two coupled nonlinear
Schr\"{o}dinger equations}, Ann. Inst. H. Poincar$\acute{e}$, Anal.
Non-Lin., 22(2005)403-439.

\bibitem{MC}
Li Ma and D.Z.Chen,  \emph{A Liouville type theorem for an integral
system} Comm. Pure and Applied Analysis, 5(2006)855-859.

\bibitem{MZ}
Li Ma, Lin Zhao, \emph{Sharp thresholds of blow up and global
existence for the coupled nonlinear Schrodinger system},Preprint,
2007.



\bibitem{SW2} E. M. Stein and G. Weiss \,\, \emph{Fractional integrals in $n$-dimensional Euclidean space} \,\,
J. Math. Mech., {\bf 7} (1958).

\bibitem{WX} J. Wei and X. Xu, \,\, \emph{Classification of solutions of
higher order conformally invariant equations}, \,\, Math. Ann.,
(1999), 207--228.

\bibitem{We1}
Weinstein, M.I., {\em Nonlinear Schr$\ddot{o}$dinger Equations and
Sharp Interpolate Estimates}, Commmu. Math. Phys., 87(1983)567-576.

\end{thebibliography}
\end{document}